\documentclass[12pt,twoside,british]{article}
\usepackage[T1]{fontenc}
\usepackage[utf8]{inputenc}
\usepackage[a4paper]{geometry}
\usepackage{amsmath}
\usepackage{amsthm}
\usepackage{amssymb}
\usepackage{setspace}

\makeatletter
\numberwithin{equation}{section}
\numberwithin{figure}{section}
\theoremstyle{plain}
\newtheorem{thm}{\protect\theoremname}[section]
\theoremstyle{plain}
\newtheorem{cor}[thm]{\protect\corollaryname}
\theoremstyle{plain}
\newtheorem{lem}[thm]{\protect\lemmaname}
\theoremstyle{definition}
\newtheorem{defn}[thm]{\protect\definitionname}
\theoremstyle{plain}
\newtheorem{prop}[thm]{\protect\propositionname}
\newcommand{\lyxaddress}[1]{
	\par {\raggedright #1
	\vspace{1.4em}
	\noindent\par}
}

\@ifundefined{date}{}{\date{}}
\usepackage[T1]{fontenc}
\usepackage[utf8]{inputenc}
\usepackage[a4paper]{geometry}
\usepackage{color}
\usepackage{dsfont}
\usepackage{amsmath}
\usepackage{amsthm}
\usepackage{amssymb}
\usepackage{setspace}
\usepackage{esint}

\makeatletter
\numberwithin{equation}{section}
\numberwithin{figure}{section}
\theoremstyle{plain}

\@ifundefined{date}{}{\date{}}
\usepackage[english]{babel}
\usepackage{latexsym}
\usepackage{amsthm}
\usepackage{esint}
\usepackage{eucal}
\usepackage{eufrak}
\usepackage{epstopdf}
\usepackage{graphicx}
\theoremstyle{plain}

\newtheoremstyle{boldremark}
    {\dimexpr\topsep/2\relax} 
    {\dimexpr\topsep/2\relax} 
    {}          
    {}          
    {\bfseries} 
    {.}         
    {.5em}      
    {}          

\theoremstyle{boldremark}
\newtheorem{brem} [thm] {Remark} 

\usepackage{fancyhdr}
\usepackage{blindtext}
\usepackage{titlesec}

\titleformat{name=\chapter}[display]
{\normalfont\Huge\itshape}
{\titlerule[1pt]\vspace{-33pt}\filleft%
  \parbox[t]{6em}{%
    \raggedleft%
    \rule{\linewidth}{0.5ex}\newline%
    \chaptertitlename\ \thechapter%
  }%
}
{4pc}
{\normalfont\upshape\bfseries\Huge}
\allowdisplaybreaks

\makeatother

\usepackage{babel}
\addto\captionsbritish{\renewcommand{\definitionname}{Definition}}
\addto\captionsbritish{\renewcommand{\lemmaname}{Lemma}}
\addto\captionsbritish{\renewcommand{\theoremname}{Theorem}}
\addto\captionsenglish{\renewcommand{\definitionname}{Definition}}
\addto\captionsenglish{\renewcommand{\lemmaname}{Lemma}}
\addto\captionsenglish{\renewcommand{\theoremname}{Theorem}}
\providecommand{\definitionname}{Definition}
\providecommand{\lemmaname}{Lemma}
\providecommand{\theoremname}{Theorem}

\makeatother

\usepackage{babel}
\providecommand{\corollaryname}{Corollary}
\providecommand{\definitionname}{Definition}
\providecommand{\lemmaname}{Lemma}
\providecommand{\propositionname}{Proposition}
\providecommand{\theoremname}{Theorem}

\begin{document}
\title{\noindent \textbf{Fractional Sobolev regularity for solutions to a
strongly degenerate parabolic equation}}
\author{\noindent Pasquale Ambrosio}
\date{November 9, 2023}
\maketitle
\begin{abstract}
\begin{singlespace}
\noindent We carry on the investigation started in \cite{Ambr2} about
the regularity of weak solutions to the strongly degenerate parabolic
equation
\[
u_{t}-\mathrm{div}\left[(\vert Du\vert-1)_{+}^{p-1}\frac{Du}{\vert Du\vert}\right]=f\,\,\,\,\,\,\,\,\,\mathrm{in}\,\,\Omega_{T}=\Omega\times(0,T),
\]
where $\Omega$ is a bounded domain in $\mathbb{R}^{n}$ for $n\geq2$,
$p\geq2$ and $\left(\,\cdot\,\right)_{+}$ stands for the positive
part. Here, we weaken the assumption on the right-hand side, by assuming
that $f\in L_{loc}^{p'}\left(0,T;B_{p',\infty,loc}^{\alpha}\left(\Omega\right)\right)$,
with $\alpha\in(0,1)$ and $p'=p/(p-1)$. This leads us to obtain
higher fractional differentiability results for a function of the
spatial gradient $Du$ of the solutions. Moreover, we establish the
higher summability of $Du$ with respect to the spatial variable.
The main novelty of the above equation is that the structure function
satisfies standard ellipticity and growth conditions only outside
the unit ball centered at the origin. We would like to point out that
the main result of this paper can be considered, on the one hand,
as the parabolic counterpart of an elliptic result contained in \cite{Ambr},
and on the other hand as the fractional version of some results established
in \cite{Ambr2}.\vspace{0.2cm}
\end{singlespace}
\end{abstract}
\begin{singlespace}
\end{singlespace}
\begin{singlespace}
\end{singlespace}
\begin{singlespace}

\section{Introduction and statement of the results}
\end{singlespace}

\begin{singlespace}
\noindent $\hspace*{1em}$In this paper, we aim to pursue our investigation
started in \cite{Ambr2} about the regularity properties of weak solutions
to the strongly degenerate parabolic equation
\begin{equation}
u_{t}-\mathrm{div}\left[(\vert Du\vert-1)_{+}^{p-1}\frac{Du}{\vert Du\vert}\right]=f\,\,\,\,\,\,\,\,\mathrm{in}\,\,\Omega_{T}=\Omega\times(0,T),\label{eq:1}
\end{equation}

\noindent where $p\geq2$, $\Omega$ is a bounded domain in $\mathbb{R}^{n}$
($n\geq2$), $\left(\,\cdot\,\right)_{+}$ stands for the positive
part and $f$ is a given function.\\
$\hspace*{1em}$The main feature of this PDE is that the structure
function satisfies standard growth and ellipticity conditions for
an exponent $p\geq2$, but only outside the unit ball centered at
the origin.\\
$\hspace*{1em}$A motivation for studying equation (\ref{eq:1}) can
be found in \textit{gas filtration problems} (for a detailed explanation,
we refer to \cite[Section 1.1]{Ambr2}).\\
$\hspace*{1em}$The elliptic version of the above equation naturally
arises in optimal transport problems with congestion effects, and
the regularity properties of its weak solutions have been widely investigated:
see, for instance, \cite{Ambr,Bog,Br0} and \cite{Br}. In this regard,
we also want to point out the very recent paper \cite{Mons}, whose
author examines the higher regularity of weak solutions to the very
degenerate elliptic system 
\[
\mathrm{div}\left[a(x)\,\frac{(\vert Dv\vert-1)_{+}^{p-1}}{\vert Dv\vert}\,Dv\right]=0\,\,\,\,\,\,\,\,\mathrm{in}\,\,\Omega,
\]
with a growth exponent $p\geq2$ and Lipschitz continuous coefficients
$a:\Omega\rightarrow\mathbb{R}$.\\
$\hspace*{1em}$To the best of our knowledge, the only parabolic counterpart
of the aforementioned works with weaker assumptions on the datum $f$
that is available in the literature is the paper \cite{Ambr2}. There,
we establish the higher differentiability of \textit{integer} order
and the higher integrability of the spatial gradient of the weak solutions
$u$ to equation (\ref{eq:1}), as well as the existence of a weak
time derivative $u_{t}$, by assuming that $f\in L^{q}\left(0,T;W^{1,q}\left(\Omega\right)\right)$
for a suitable exponent $q>1$.\\
$\hspace*{1em}$Here, instead, we prove local higher differentiability
results for a certain function of the spatial gradient $Du$ in the
scale of \textit{fractional} Sobolev spaces, under the assumption
that the datum $f$ belongs to the local Bochner space $L_{loc}^{p'}\left(0,T;B_{p',\infty,loc}^{\alpha}\left(\Omega\right)\right)$,
where $p'=p/(p-1)$ is the conjugate exponent of $p$, $\alpha\in(0,1)$,
while $B_{p',\infty}^{\alpha}\left(\Omega\right)$ denotes a particular
class of Besov functions (see Section \ref{sec:Besov} below for the
definition).\\
Furthermore, we establish the local higher summability of $Du$ with
respect to the spatial variables, under the same assumption on the
regularity of $f$.\\
$\hspace*{1em}$The distinguishing feature of equation (\ref{eq:1})
is that the principal part behaves like a \textit{$p$-Laplace operator
only at infinity}. Before giving the main result, let us summarize
a few previous results on this topic: the regularity of solutions
to parabolic problems with asymptotic structure of $p$-Laplacian
type has been explored in \cite{Isernia}, where a BMO regularity
has been proved for solutions to asymptotically parabolic systems
in the case $p=2$ and $f=0$ (see also \cite{Kuusi}, where the local
Lipschitz continuity of weak solutions with respect to the spatial
variable is established). In addition, we want to mention the results
contained in \cite{Byun}, where nonhomogeneous parabolic problems
involving a discontinuous nonlinearity and an asymptotic regularity
over an irregular domain in divergence form of $p$-Laplacian type
are considered. There, the authors establish a global Calderón-Zygmund
estimate by converting a given asymptotically regular problem to a
suitable regular problem.\\
$\hspace*{1em}$We would like to notice that our assumption on the
datum $f$ is weaker than those considered in the mentioned papers.
This prevents us from achieving a (local) higher summability for $Du$
over $\Omega_{T}$ (i.e. with respect to the space-time variable $z=(x,t)$),
contrary to what we were able to do in \cite[Theorem 1.3]{Ambr2}.
That is why here we obtain the higher summability of $Du$ only with
respect to the spatial variable (see Corollary \ref{cor:cor1} below).
For the same reason, we are led to deal with fractional Sobolev spaces,
rather than the traditional Sobolev spaces of integer order.\\
$\hspace*{1em}$The main result we prove in this paper is the following
theorem. We refer to Sections \ref{sec:prelim} and \ref{sec:Functional spaces}
for notation and definitions.\\

\end{singlespace}
\begin{thm}
\begin{singlespace}
\noindent \label{thm:mainth} Let $n\geq2$, $p\geq2$, $\alpha\in(0,1)$
and $f\in L_{loc}^{p'}\left(0,T;B_{p',\infty,loc}^{\alpha}\left(\Omega\right)\right)$.
Moreover, assume that 
\[
u\in C^{0}\left(0,T;L^{2}\left(\Omega\right)\right)\cap L^{p}\left(0,T;W^{1,p}\left(\Omega\right)\right)
\]
is a weak solution of equation $\mathrm{(\ref{eq:1})}$. Then the
solution satisfies 
\[
H_{\frac{p}{2}}(Du)\,\in\,L_{loc}^{2}\left(0,T;W_{loc}^{s,2}\left(\Omega,\mathbb{R}^{n}\right)\right)\,\,\,\,\,\,\,\,\,for\,\,all\,\,s\in\left(0,\frac{\alpha+1}{2}\right),
\]

\noindent where 
\begin{equation}
H_{\frac{p}{2}}(Du):=\left(\left|Du\right|-1\right)_{+}^{p/2}\frac{Du}{\left|Du\right|}.\label{eq:H-fun}
\end{equation}
Furthermore, the following local estimate \begin{align*}
&\sup_{t_{0}-(\varrho/2)^{2}<t<t_{0}}\int_{B_{\varrho/2}(x_{0})}\left|\tau_{i,h}u(x,t)\right|^{2}dx\,+\,\int_{Q_{\varrho/2}(z_{0})}\left|\tau_{i,h}H_{\frac{p}{2}}(Du)\right|^{2}dz\\
&\,\,\,\,\,\,\,\leq\,c\,\left[\left|h\right|^{2}\int_{Q_{2\varrho}(z_{0})}\left(\left|Du\right|^{p}+1\right)\,dz\,+\,\left|h\right|^{\alpha+1}\Vert Du\Vert_{L^{p}(Q_{2\varrho})}\,\Vert f\Vert_{L^{p'}\left(t_{0}\,-\,4\varrho^{2},\,t_{0}\,;\,B_{p',\infty}^{\alpha}\left(B_{2\varrho}\right)\right)}\right]
\end{align*}holds true for any $i\in\left\{ 1,\ldots,n\right\} $, for any parabolic
cylinder $Q_{\varrho}(z_{0})\subset Q_{2\varrho}(z_{0})\Subset\Omega_{T}$,
for any $h\in\mathbb{R}$ such that $\left|h\right|<\varrho/4$ and
a positive constant $c$ depending at most on $p$, $n$ and $\varrho$. 
\end{singlespace}
\end{thm}

\begin{singlespace}
\noindent $\hspace*{1em}$The proof of Theorem \ref{thm:mainth} is
achieved using the well-known difference quotients technique in the
spatial directions (see Section \ref{sec:proofs} below). Here we
will argue as in \cite[Lemma 5.1]{Duzaar} and \cite[Theorem 4.1]{Giann},
but we need to take into account the strong degeneracy of equation
(\ref{eq:1}), exactly as we have done in \cite{Ambr2}. This is why
we obtain the higher fractional differentiability not for the spatial
gradient of the solution itself, but for a function of the spatial
gradient $Du$ that vanishes in the set where equation (\ref{eq:1})
becomes degenerate.\\
$\hspace*{1em}$Actually, Theorem \ref{thm:mainth} can be considered
as the parabolic counterpart of Theorem 4.1 in \cite{Ambr}, where
however we obtained a higher fractional differentiability result in
the scale of \textit{Besov} spaces, starting from a datum in the local
Besov-Lipschitz space $B_{p',\infty,loc}^{\alpha}\left(\Omega\right)$,
for some $\alpha\in(0,1)$. Furthermore, the above theorem can also
be viewed as the fractional version of Theorems 1.1 and 1.4 in \cite{Ambr2}.\\
$\hspace*{1em}$As a consequence of the previous result, using a Sobolev-type
embedding theorem for fractional Sobolev spaces, we get a gain of
summability for $Du$ with respect to the spatial variable. More precisely,
we have the following 
\end{singlespace}
\begin{cor}
\begin{singlespace}
\noindent \label{cor:cor1} Under the assumptions of Theorem \ref{thm:mainth},
we obtain that
\[
H_{\frac{p}{2}}(Du)\,\in\,L_{loc}^{2}\left(0,T;L_{loc}^{r}\left(\Omega,\mathbb{R}^{n}\right)\right)\,\,\,\,\,\,\,\,\,for\,\,all\,\,r\in\left[1,\frac{2n}{n-\alpha-1}\right)
\]
and
\[
Du\,\in\,L_{loc}^{p}\left(0,T;L_{loc}^{q}\left(\Omega,\mathbb{R}^{n}\right)\right)\,\,\,\,\,\,\,\,\,for\,\,all\,\,q\in\left[1,\frac{np}{n-\alpha-1}\right).
\]
\end{singlespace}
\end{cor}

\begin{singlespace}
\noindent $\hspace*{1em}$It is worth pointing out that, starting
from the assumption $f\in L_{loc}^{p'}\left(0,T;B_{p',\infty,loc}^{\alpha}\left(\Omega\right)\right)$,
higher fractional differentiability results such as that of Theorem
\ref{thm:mainth} above seem not to have been established yet for
the spatial gradient of weak solutions to the strongly degenerate
equation (\ref{eq:1}).\\
$\hspace*{1em}$However, we would like to mention the result contained
in \cite[Theorem 2.10]{Byun}, whose authors consider a weak solution
$u\in C^{0}\left([0,T];L^{2}\left(\Omega\right)\right)\cap L^{p}\left(0,T;W_{0}^{1,p}\left(\Omega\right)\right)$
to the following nonlinear parabolic problem in divergence form:
\[
\begin{cases}
\begin{array}{cc}
u_{t}-\mathrm{div}\,\mathbf{a}(x,t,Du)=\mathrm{div}\left(\left|F\right|^{p-2}F\right) & \,\,\,\mathrm{in}\,\,\,\Omega_{T},\\
u=0 & \,\,\,\,\,\,\,\,\,\,\,\,\,\,\mathrm{on}\,\,\,\partial_{\mathrm{par}}\Omega_{T},
\end{array}\end{cases}
\]
where $\frac{2n}{n+2}<p<\infty$, $\mathbf{a}$ is a discontinuous
nonlinearity with an asymptotic regularity, $\Omega$ is a bounded
domain whose boundary $\partial\Omega$ is nonsmooth and $F=F(x,t)=\left(f_{1}(x,t),\ldots,f_{n}(x,t)\right)\in L^{p}(\Omega_{T},\mathbb{R}^{n})$
is a given vector-valued function. In fact, the just mentioned result
gives an affirmative answer as to what are both the weakest regularity
requirement on $\mathbf{a}$ and the lowest level of geometric assumption
on $\partial\Omega$ under which the implication 
\[
\left|F\right|^{p}\in L^{q}(\Omega_{T})\Longrightarrow\left|Du\right|^{p}\in L^{q}(\Omega_{T})
\]
holds true for every $q\in[1,\infty)$.
\end{singlespace}
\begin{singlespace}

\subsection{Comparison with a less degenerate elliptic equation}
\end{singlespace}

\begin{singlespace}
\noindent $\hspace*{1em}$Before describing the structure of this
paper, we wish to make some considerations about the interplay between
the regularity of the right-hand side of (\ref{eq:1}) and that of
the vector field $H_{\frac{p}{2}}(Du)$ defined in (\ref{eq:H-fun}),
starting from the comparison with a known result for an elliptic equation
which is less degenerate than (\ref{eq:1}).\\
$\hspace*{1em}$In \cite[Theorem 1.1]{BraSan}, Brasco and Santambrogio
established the sharp assumptions on the datum $g$ in order to obtain
that the $W^{1,p}$ solutions to the Poisson equation for the $p$-Laplace
operator
\begin{equation}
-\Delta_{p}v:=-\,\mathrm{div}(\left|Dv\right|^{p-2}Dv)=g\,\,\,\,\,\,\,\,\mathrm{in}\,\,\Omega,\label{eq:p-Laplace}
\end{equation}
still satisfy the following Uhlenbeck's result (see \cite[Lemma 3.1]{Uhl})
\begin{equation}
\left|Dv\right|^{\frac{p-2}{2}}Dv\in W_{loc}^{1,2}(\Omega,\mathbb{R}^{n}),\label{eq:Uhlenbeck}
\end{equation}
in the superquadratic case $p>2$, where $\Omega\subset\mathbb{R}^{n}$
is an open set. For the sake of completeness, we recall their result
below here, by noting that $g$ must belong at least to a suitable
(local) fractional Sobolev space for condition (\ref{eq:Uhlenbeck})
to be true:\medskip{}

\end{singlespace}
\begin{thm}
\begin{singlespace}
\noindent \textbf{\textup{(\cite[Theorem 1.1]{BraSan}).}} \label{thm:BraSan}
Let $p>2$ and let $U\in W_{loc}^{1,p}(\Omega)$ be a local weak solution
of equation $(\ref{eq:p-Laplace})$. If 
\begin{equation}
g\in W_{loc}^{s,p'}(\Omega)\,\,\,\,\,\,\,\,\,\,\,\,\,with\,\,\,\,\,\frac{p-2}{p}<s\leq1,\label{eq:sharp assump}
\end{equation}
then 
\[
\mathcal{V}:=\left|DU\right|^{\frac{p-2}{2}}DU\in W_{loc}^{1,2}(\Omega,\mathbb{R}^{n}),
\]
and
\[
DU\in W_{loc}^{\sigma,p}(\Omega,\mathbb{R}^{n}),\,\,\,\,\,\,\,\,\,\,\,\,\,for\,\,\,\,\,0<\sigma<\frac{2}{p}.
\]
\end{singlespace}
\end{thm}

\begin{singlespace}
\noindent As for the sharpness of assumption (\ref{eq:sharp assump}),
we refer the interested reader to \cite[Sections 1.2 and 5]{BraSan}.
Under the hypotheses of Theorem \ref{thm:BraSan}, we find that $DU$
also belongs to some fractional Sobolev space (locally in $\Omega$).\\
$\hspace*{1em}$The interplay between the regularity of the right-hand
side $g$ and that of the vector field $\mathcal{V}$ has been considered
in detail also in \cite{Ming2}, where the point of view is slightly
different: the main concern there is to obtain (fractional) differentiability
of the vector field $\mathcal{V}$ when $g$ is not regular. In particular,
in \cite{Ming2} the datum $g$ may not belong to the relevant dual
Sobolev space and the concept of solution to (\ref{eq:p-Laplace})
has to be carefully defined. We refer to \cite[Remark 1.3]{BraSan}
and the references therein for comparison with other results.\\
$\hspace*{1em}$In the case of the strongly degenerate parabolic equation
(\ref{eq:1}), the vector-valued function $\mathcal{V}$ is replaced
by the function $H_{\frac{p}{2}}(Du)$ defined in (\ref{eq:H-fun}),
which belongs to $L_{loc}^{2}\left(0,T;W_{loc}^{s,2}\left(\Omega,\mathbb{R}^{n}\right)\right)$
for all $s\in\left(0,\frac{\alpha+1}{2}\right)$ under the assumptions
of the Theorem \ref{thm:mainth} that we prove in this paper. One
could then ask himself what are the optimal (weakest) assumptions
to be imposed on the datum $f$ so that the weak solutions of equation
(\ref{eq:1}) still satisfy ``our'' following condition (see \cite[Theorems 1.1 and 1.4]{Ambr2})
\begin{equation}
H_{\frac{p}{2}}(Du)\,\in\,L_{loc}^{2}(0,T;W_{loc}^{1,2}\left(\Omega,\mathbb{R}^{n}\right)),\label{eq:previous}
\end{equation}
for $p\geq2$, since here we cannot achieve the result (\ref{eq:previous})
under the hypothesis
\[
f\in L_{loc}^{p'}\left(0,T;B_{p',\infty,loc}^{\alpha}\left(\Omega\right)\right),\,\,\,\,\,\,\,\,\,\,\mathrm{with}\,\,\,\,\,\alpha\in(0,1).
\]
$\hspace*{1em}$However, it is worth emphasizing that equation (\ref{eq:1})
exhibits a more severe degeneracy than equation (\ref{eq:p-Laplace})
and that the crucial point in \cite{BraSan} that 
\[
\left|DU\right|^{\frac{p-2}{2}}DU\in W_{loc}^{1,2}\,\,\Longrightarrow\,\,DU\in W_{loc}^{\sigma,p}
\]
cannot be retrieved in our framework.\\
\\
$\hspace*{1em}$The paper is organized as follows. Section \ref{sec:prelim}
is devoted to the preliminaries: after a list of some classic notations
and some essentials estimates, we recall the fundamental properties
of the difference quotients of Sobolev functions. In Section \ref{sec:Functional spaces}
we recall the basic facts on the functional spaces involved in this
paper: subsection \ref{sec:fracSob} is entirely devoted to the definition
and properties of the fractional Sobolev spaces that will be useful
to establish our results, while in subsection \ref{sec:Besov} we
give the definition of the Besov spaces $B_{q,\infty}^{\alpha}\left(\Omega\right)$
for $0<\alpha<1$ and $1\leq q<\infty$, in order to introduce the
Bochner space $L^{q}\left(0,T;B_{q,\infty}^{\alpha}\left(\Omega\right)\right)$
and its local version. Finally, in Section \ref{sec:proofs} we prove
Theorem \ref{thm:mainth} and Corollary \ref{cor:cor1}.
\end{singlespace}
\begin{singlespace}

\section{Notations and preliminaries\label{sec:prelim}}
\end{singlespace}

\begin{singlespace}
\noindent $\hspace*{1em}$In this paper we shall denote by $C$ or
$c$ a general positive constant that may vary on different occasions.
Relevant dependencies on parameters and special constants will be
suitably emphasized using parentheses or subscripts. The norm we use
on $\mathbb{R}^{n}$ will be the standard Euclidean one and it will
be denoted by $\left|\,\cdot\,\right|$. In particular, for the vectors
$\xi,\eta\in\mathbb{R}^{n}$, we write $\langle\xi,\eta\rangle$ for
the usual inner product and $\left|\xi\right|:=\langle\xi,\xi\rangle^{\frac{1}{2}}$
for the corresponding Euclidean norm.\\
$\hspace*{1em}$For points in space-time, we will frequently use abbreviations
like $z=(x,t)$ or $z_{0}=(x_{0},t_{0})$, for spatial variables $x$,
$x_{0}\in\mathbb{R}^{n}$ and times $t$, $t_{0}\in\mathbb{R}$. We
also denote by $B(x_{0},\varrho)=B_{\varrho}(x_{0})=\left\{ x\in\mathbb{R}^{n}:\left|x-x_{0}\right|<\varrho\right\} $
the open ball with radius $\varrho>0$ and center $x_{0}\in\mathbb{R}^{n}$;
when not important, or clear from the context, we shall omit to denote
the center as follows: $B_{\varrho}\equiv B(x_{0},\varrho)$. Unless
otherwise stated, different balls in the same context will have the
same center. Moreover, we use the notation
\[
Q_{\varrho}(z_{0}):=B_{\varrho}(x_{0})\times(t_{0}-\varrho^{2},t_{0}),\,\,\,\,\,z_{0}=(x_{0},t_{0})\in\mathbb{R}^{n}\times\mathbb{R},\,\,\varrho>0,
\]
for the backward parabolic cylinder with vertex $(x_{0},t_{0})$ and
width $\varrho$. We shall sometimes omit the dependence on the vertex
when all cylinders occurring in a proof share the same vertex.\\
\\
$\hspace*{1em}$We now recall some tools that will be useful to prove
our results. For the auxiliary function $H_{\lambda}:\mathbb{R}^{n}\rightarrow\mathbb{R}^{n}$
defined as
\[
H_{\lambda}(\xi):=\begin{cases}
\begin{array}{cc}
\left(\left|\xi\right|-1\right)_{+}^{\lambda}\frac{\xi}{\left|\xi\right|} & \,\,\,\,\,\,\,\,\,\,\,\,\,\,\,\,\,\,\,\mathrm{if}\,\,\,\xi\in\mathbb{R}^{n}\setminus\left\{ 0\right\} ,\\
0 & \mathrm{if}\,\,\,\xi=0,
\end{array}\end{cases}
\]
where $\lambda>0$ is a parameter, we record the following estimates
(see \cite[Lemma 4.1]{Br}):\\

\end{singlespace}
\begin{lem}
\begin{singlespace}
\noindent \label{lem:Brasco} If $2\leq p<\infty$, then for every
$\xi,\eta\in\mathbb{R}^{n}$ we get 
\[
\langle H_{p-1}(\xi)-H_{p-1}(\eta),\xi-\eta\rangle\,\geq\,\frac{4}{p^{2}}\left|H_{\frac{p}{2}}(\xi)-H_{\frac{p}{2}}(\eta)\right|^{2},
\]
\[
\left|H_{p-1}(\xi)-H_{p-1}(\eta)\right|\,\leq\,(p-1)\left(\left|H_{\frac{p}{2}}(\xi)\right|^{\frac{p-2}{p}}+\left|H_{\frac{p}{2}}(\eta)\right|^{\frac{p-2}{p}}\right)\left|H_{\frac{p}{2}}(\xi)-H_{\frac{p}{2}}(\eta)\right|.
\]
\end{singlespace}
\end{lem}

\begin{singlespace}
\noindent $\hspace*{1em}$

\noindent $\hspace*{1em}$We conclude this first part of the preliminaries
by recalling the following\vspace{0.04mm}

\end{singlespace}
\begin{defn}
\begin{singlespace}
\noindent A function $u\in C^{0}\left(0,T;L^{2}\left(\Omega\right)\right)\cap L^{p}\left(0,T;W^{1,p}\left(\Omega\right)\right)$
is a \textit{weak solution} of equation (\ref{eq:1}) iff for any
test function $\varphi\in C_{0}^{\infty}(\Omega_{T})$ the following
integral identity holds:
\begin{equation}
\int_{\Omega_{T}}\left(u\cdot\partial_{t}\varphi-\langle H_{p-1}(Du),D\varphi\rangle\right)\,dz\,=\,-\int_{\Omega_{T}}f\varphi\,dz.\label{eq:weaksol}
\end{equation}
\end{singlespace}
\end{defn}

\begin{singlespace}

\subsection{Difference quotients\label{subsec:DiffOpe}}
\end{singlespace}

\begin{singlespace}
\noindent $\hspace*{1em}$We recall here the definition and some elementary
properties of the difference quotients that will be useful in the
following (see, for instance, \cite{Giu}).\\

\end{singlespace}
\begin{defn}
\begin{singlespace}
\noindent For every vector-valued function $F:\mathbb{R}^{n}\rightarrow\mathbb{R}^{N}$
the \textit{finite difference operator }in the direction $x_{i}$
is defined by
\[
\tau_{i,h}F(x)=F(x+he_{i})-F(x),
\]
where $h\in\mathbb{R}$, $e_{i}$ is the unit vector in the direction
$x_{i}$ and $i\in\left\{ 1,\ldots,n\right\} $. \\
$\hspace*{1em}$The \textit{difference quotient} of $F$ with respect
to $x_{i}$ is defined for $h\in\mathbb{R}\setminus\left\{ 0\right\} $
by 
\[
\Delta_{i,h}F(x)=\frac{\tau_{i,h}F(x)}{h}.
\]
\end{singlespace}
\end{defn}

\begin{singlespace}
\noindent $\hspace*{1em}$When no confusion can arise, we shall omit
the index $i$ and simply write $\tau_{h}$ or $\Delta_{h}$ instead
of $\tau_{i,h}$ or $\Delta_{i,h}$, respectively. 
\end{singlespace}
\begin{prop}
\begin{singlespace}
\noindent Let $F$ be a function such that $F\in W^{1,q}\left(\Omega\right)$,
with $q\geq1$, and let us consider the set
\[
\Omega_{\left|h\right|}:=\left\{ x\in\Omega:\mathrm{dist}\left(x,\partial\Omega\right)>\left|h\right|\right\} .
\]
Then:\\
\\
$\mathrm{(}a\mathrm{)}$ $\Delta_{h}F\in W^{1,q}\left(\Omega_{\left|h\right|}\right)$
and
\[
D_{j}(\Delta_{h}F)=\Delta_{h}(D_{j}F),\,\,\,\,\,for\,\,every\,\,j\in\left\{ 1,\ldots,n\right\} .
\]
$\mathrm{(}b\mathrm{)}$ If at least one of the functions $F$ or
$G$ has support contained in $\Omega_{\left|h\right|}$, then
\[
\int_{\Omega}F\,\Delta_{h}G\,dx\,=\,-\int_{\Omega}G\,\Delta_{-h}F\,dx.
\]
$\mathrm{(}c\mathrm{)}$ We have 
\[
\Delta_{h}(FG)(x)=F(x+he_{i})\Delta_{h}G(x)\,+\,G(x)\Delta_{h}F(x).
\]
\vspace{0.01mm}
\end{singlespace}
\end{prop}

\begin{singlespace}
\noindent $\hspace*{1em}$The next result about the finite difference
operator is a kind of integral version of the Lagrange Theorem.
\end{singlespace}
\begin{lem}
\begin{singlespace}
\noindent \label{lem:Giusti1} If $0<\varrho<R$, $\left|h\right|<\frac{R-\varrho}{2}$,
$1<q<+\infty$, and $F\in L^{q}\left(B_{R},\mathbb{R}^{N}\right)$,
$DF\in L^{q}\left(B_{R},\mathbb{R}^{N\times n}\right)$, then
\[
\int_{B_{\varrho}}\left|\tau_{h}F(x)\right|^{q}dx\,\leq\,c^{q}(n)\left|h\right|^{q}\int_{B_{R}}\left|DF(x)\right|^{q}dx.
\]
Moreover
\[
\int_{B_{\varrho}}\left|F(x+he_{i})\right|^{q}dx\,\leq\,\int_{B_{R}}\left|F(x)\right|^{q}dx.
\]
\end{singlespace}
\end{lem}

\begin{singlespace}
\noindent $\hspace*{1em}$Finally, we recall the following fundamental
result, whose proof can be found in \cite[Lemma 8.2]{Giu}:
\end{singlespace}
\begin{lem}
\begin{singlespace}
\noindent \label{lem:RappIncre} Let $F:\mathbb{R}^{n}\rightarrow\mathbb{R}^{N}$,
$F\in L^{q}\left(B_{R},\mathbb{R}^{N}\right)$ with $1<q<+\infty$.
Suppose that there exist $\varrho\in(0,R)$ and a constant $M>0$
such that 
\[
\sum_{i=1}^{n}\int_{B_{\varrho}}\left|\tau_{i,h}F(x)\right|^{q}dx\,\leq\,M^{q}\left|h\right|^{q}
\]
for every $h$ with $\left|h\right|<\frac{R-\varrho}{2}$. Then $F\in W^{1,q}\left(B_{\varrho},\mathbb{R}^{N}\right)$.
Moreover 
\[
\Vert DF\Vert_{L^{q}\left(B_{\varrho}\right)}\leq M
\]
and
\[
\Delta_{i,h}F\rightarrow D_{i}F\,\,\,\,\,\,\,\,\,\,in\,\,L_{loc}^{q}\left(B_{R}\right),\,\,as\,\,h\rightarrow0,
\]
for each $i\in\left\{ 1,\ldots,n\right\} $.
\end{singlespace}
\end{lem}

\begin{singlespace}
\noindent 

\end{singlespace}\begin{singlespace}

\section{Functional spaces \label{sec:Functional spaces}}
\end{singlespace}

\begin{singlespace}
\noindent $\hspace*{1em}$Here we recall some essential facts about
the functional spaces involved in this paper, starting with the definition
and some properties of the fractional Sobolev spaces that will be
useful to prove our results (see, for instance, \cite{DiNezza}).
\end{singlespace}
\begin{singlespace}

\subsection{Fractional Sobolev spaces \label{sec:fracSob}}
\end{singlespace}

\begin{singlespace}
\noindent $\hspace*{1em}$Let $\Omega$ be a general, possibly nonsmooth,
bounded open set in $\mathbb{R}^{n}$. For any $s\in(0,1)$ and for
any $q\in[1,+\infty)$, we define the fractional Sobolev space $W^{s,q}\left(\Omega,\mathbb{R}^{k}\right)$
as follows
\[
W^{s,q}\left(\Omega,\mathbb{R}^{k}\right):=\left\{ v\in L^{q}\left(\Omega,\mathbb{R}^{k}\right):\frac{\left|v(x)-v(y)\right|}{\left|x-y\right|^{\frac{n}{q}\,+\,s}}\in L^{q}\left(\Omega\times\Omega\right)\right\} ,
\]
i.e. an intermediate Banach space between $L^{q}\left(\Omega,\mathbb{R}^{k}\right)$
and $W^{1,q}\left(\Omega,\mathbb{R}^{k}\right)$, endowed with the
norm 
\[
\Vert v\Vert_{W^{s,q}(\Omega)}:=\left(\int_{\Omega}\left|v\right|^{q}dx\,+\,\int_{\Omega}\int_{\Omega}\frac{\left|v(x)-v(y)\right|^{q}}{\left|x-y\right|^{n\,+\,sq}}\,dx\,dy\right)^{\frac{1}{q}},
\]
where the term
\[
\left[v\right]_{W^{s,q}(\Omega)}:=\left(\int_{\Omega}\int_{\Omega}\frac{\left|v(x)-v(y)\right|^{q}}{\left|x-y\right|^{n\,+\,sq}}\,dx\,dy\right)^{\frac{1}{q}}
\]
is the so-called \textit{Gagliardo seminorm} of $v$.\\
$\hspace*{1em}$As in the classic case with $s$ being an integer,
the space $W^{s',q}\left(\Omega\right)$ is continuously embedded
in $W^{s,q}\left(\Omega\right)$ when $s\leq s'$, as shown by the
next result (see \cite[Proposition 2.1]{DiNezza}).\medskip{}

\end{singlespace}
\begin{prop}
\begin{singlespace}
\noindent Let $\Omega$ be an open set in $\mathbb{R}^{n}$, $q\in[1,+\infty)$
and $0<s\leq s'<1$. Then there exists a constant $C\equiv C(n,s,q)\geq1$
such that, for any $v\in W^{s',q}\left(\Omega\right)$, we have 
\[
\Vert v\Vert_{W^{s,q}(\Omega)}\leq\,C\,\Vert v\Vert_{W^{s',q}(\Omega)}.
\]
In particular, $W^{s',q}\left(\Omega\right)\subseteq W^{s,q}\left(\Omega\right)$. 
\end{singlespace}
\end{prop}

\begin{singlespace}
\noindent $\hspace*{1em}$As is well known when $s\in\mathbb{N}$,
under certain regularity assumptions on the open set $\Omega\subset\mathbb{R}^{n}$,
any function in $W^{s,q}\left(\Omega\right)$ can be extended to a
function in $W^{s,q}\left(\mathbb{R}^{n}\right)$. Extension results
are needed to improve some embedding theorems, in the classic case
as well as in the fractional one. In this regard, we now give the
following
\end{singlespace}
\begin{defn}
\begin{singlespace}
\noindent For any $s\in(0,1)$ and any $q\in[1,\infty)$, we say that
an open set $\Omega\subseteq\mathbb{R}^{n}$ is an \textit{extension
domain} for $W^{s,q}$ if there exists a positive constant $C\equiv C(n,q,s,\Omega)$
such that: for every function $v\in W^{s,q}\left(\Omega\right)$ there
exists $\widetilde{v}\in W^{s,q}\left(\mathbb{R}^{n}\right)$ with
$\widetilde{v}|_{\Omega}=v$ and
\[
\Vert\widetilde{v}\Vert_{W^{s,q}(\mathbb{R}^{n})}\leq\,C\,\Vert v\Vert_{W^{s,q}(\Omega)}.
\]
\end{singlespace}
\end{defn}

\begin{singlespace}
\noindent $\hspace*{1em}$In general, an arbitrary open set may not
be an extension domain for $W^{s,q}$. However, the following result
ensures that every open Lipschitz set $\Omega$ with bounded boundary
is an extension domain for $W^{s,q}$ (a proof can be found in \cite[Theorem 5.4]{DiNezza}). 
\end{singlespace}
\begin{thm}
\begin{singlespace}
\noindent \label{thm:Lipschitz} Let $q\in[1,+\infty)$, $s\in(0,1)$
and $\Omega\subseteq\mathbb{R}^{n}$ be an open set of class $C^{0,1}$
with bounded boundary. Then $W^{s,q}\left(\Omega\right)$ is continuously
embedded in $W^{s,q}\left(\mathbb{R}^{n}\right)$, namely for any
$v\in W^{s,q}\left(\Omega\right)$ there exists $\widetilde{v}\in W^{s,q}\left(\mathbb{R}^{n}\right)$
such that $\widetilde{v}|_{\Omega}=v$ and 
\[
\Vert\widetilde{v}\Vert_{W^{s,q}(\mathbb{R}^{n})}\leq\,C\,\Vert v\Vert_{W^{s,q}(\Omega)}
\]
for some positive constant $C\equiv C(n,q,s,\Omega)$. 
\end{singlespace}
\end{thm}

\begin{singlespace}
\noindent For more information on the problem of characterizing the
class of sets that are extension domains for $W^{s,q}$, we refer
the interested reader to Zhou's paper \cite{Zhou}, where an answer
to this question has been given (see also \cite{Jones} and \cite[Chapters 11 and 12]{Leoni}).\\
$\hspace*{1em}$For further needs, we now recall the following Sobolev-type
embedding theorem, whose proof can be found in \cite[Theorem 6.7]{DiNezza}.
\end{singlespace}
\begin{thm}
\begin{singlespace}
\noindent \label{thm:fracemb} Let $s\in(0,1)$ and $q\in[1,+\infty)$
be such that $sq<n$. Let $\Omega\subseteq\mathbb{R}^{n}$ be an extension
domain for $W^{s,q}$. Then there exists a positive constant $C\equiv C(n,q,s,\Omega)$
such that, for any $v\in W^{s,q}\left(\Omega\right)$, we have
\[
\Vert v\Vert_{L^{r}(\Omega)}\leq\,C\,\Vert v\Vert_{W^{s,q}(\Omega)}
\]
for any $r\in\left[q,q^{*}\right]$; i.e. the space $W^{s,q}\left(\Omega\right)$
is continuously embedded in $L^{r}\left(\Omega\right)$ for any $r\in\left[q,q^{*}\right]$,
where $q^{*}:=nq/(n-sq)$ is the so-called ``fractional critical
exponent''.\\
$\hspace*{1em}$Moreover, if $\Omega$ is bounded, then the space
$W^{s,q}\left(\Omega\right)$ is continuously embedded in $L^{r}\left(\Omega\right)$
for any $r\in\left[1,q^{*}\right]$.
\end{singlespace}
\end{thm}

\begin{singlespace}
\noindent \begin{brem} In the critical case $r=q^{*}$, the constant
$C$ in Theorem \ref{thm:fracemb} does not depend on $\Omega$ (see
Remark 6.8 in \cite{DiNezza}).\end{brem}\medskip{}

\noindent $\hspace*{1em}$For the treatment of parabolic equations,
the following function space plays an important role. Let $1\leq q<\infty$
and $0<s<1$. A map $g\in L^{q}\left(\Omega\times(t_{0},t_{1}),\mathbb{R}^{k}\right)$
belongs to the space $L^{q}\left(t_{0},t_{1};W^{s,q}\left(\Omega,\mathbb{R}^{k}\right)\right)$
if and only if
\[
\int_{t_{0}}^{t_{1}}\int_{\Omega}\int_{\Omega}\frac{\left|g(x,t)-g(y,t)\right|^{q}}{\left|x-y\right|^{n\,+\,sq}}\,dx\,dy\,dt<\infty.
\]
In this paper, we will use the corresponding local version of this
space, which will be denoted by the subscript ``\textit{loc}''.
More precisely, we write $g\in L_{loc}^{q}\left(0,T;W_{loc}^{s,q}\left(\Omega,\mathbb{R}^{k}\right)\right)$
if and only if $g\in L^{q}\left(t_{0},t_{1};W^{s,q}\left(\Omega',\mathbb{R}^{k}\right)\right)$
for all domains $\Omega'\times(t_{0},t_{1})\Subset\Omega_{T}$.\\
$\hspace*{1em}$Now we conclude this section with the parabolic version
of the well-known result about the relation between Nikolskii spaces
and fractional Sobolev spaces. This result is contained in \cite[Lemma 2.4]{Sche}
(see also \cite[Proposition 2.19]{Duzaar}), and its proof can be
obtained by a simple adaptation of the standard elliptic results \cite{DuGaMi,EspLeoMin,Ming1,Ming2,Ming3}. 
\end{singlespace}
\begin{prop}
\begin{singlespace}
\noindent \label{prop:Nikolskii} Let $Q_{\sigma}(z_{0})\subset\mathbb{R}^{n+1}$
be a parabolic cylinder. Moreover, assume that $G\in L^{q}\left(Q_{\sigma}(z_{0}),\mathbb{R}^{k}\right)$,
where $1\leq q<\infty$. Then, for any $\theta\in(0,1)$, the estimate
\[
\left|h\right|^{-q\theta}\int_{Q_{\sigma}(z_{0})}\left|G(x+he_{i},t)-G(x,t)\right|^{q}dx\,dt\,\leq\,M^{q}\,<\infty
\]
for a fixed constant $M\geq0$, every $0\neq\left|h\right|\leq h_{0}$
and every $i\in\left\{ 1,\ldots,n\right\} $ implies 
\[
G\,\in\,L_{loc}^{q}\left(t_{0}-\sigma^{2},t_{0};W_{loc}^{s,q}\left(B_{\sigma}(x_{0}),\mathbb{R}^{k}\right)\right)
\]
for all $s\in(0,\theta)$.
\end{singlespace}
\end{prop}

\begin{singlespace}

\subsection{Besov spaces \label{sec:Besov}}
\end{singlespace}

\begin{singlespace}
\noindent $\hspace*{1em}$Here we recall the definition of the Besov
space $B_{q,\infty}^{\alpha}\left(\Omega\right)$ for $0<\alpha<1$
and $1\leq q<\infty$ (see Section 2.5.12 in \cite{Tri}). For a function
$v\in L^{q}\left(\Omega\right)$ we say that $v\in B_{q,\infty}^{\alpha}\left(\Omega\right)$
if 
\begin{equation}
\left[v\right]_{\dot{B}_{q,\infty}^{\alpha}\left(\Omega\right)}:=\sup_{h\in\mathbb{R}^{n}}\left(\int_{\Omega}\frac{\left|\Delta_{h}[v](x)\right|^{q}}{\left|h\right|^{\alpha q}}\,dx\right)^{\frac{1}{q}}<\infty,\label{eq:seminorm}
\end{equation}
where $\Delta_{h}[v](x):=[v(x+h)-v(x)]\cdot\mathbf{{1}}_{\Omega}(x+h)$.
One can define a norm on the space $B_{q,\infty}^{\alpha}\left(\Omega\right)$
as follows
\[
\Vert v\Vert_{B_{q,\infty}^{\alpha}\left(\Omega\right)}:=\left(\Vert v\Vert_{L^{q}\left(\Omega\right)}^{q}+\left[v\right]_{\dot{B}_{q,\infty}^{\alpha}\left(\Omega\right)}^{q}\right)^{\frac{1}{q}},
\]
and with this norm $B_{q,\infty}^{\alpha}\left(\Omega\right)$ is
a Banach space. Actually, in (\ref{eq:seminorm}) we can simply take
the supremum over $\left|h\right|<\delta$ for a fixed $\delta>0$
and obtain an equivalent norm, since 
\[
\sup_{\left|h\right|\geq\delta}\left(\int_{\Omega}\frac{\left|\Delta_{h}[v](x)\right|^{q}}{\left|h\right|^{\alpha q}}\,dx\right)^{\frac{1}{q}}\leq c(n,\alpha,q,\delta)\,\Vert v\Vert_{L^{q}\left(\Omega\right)}.
\]
By construction, we have $B_{q,\infty}^{\alpha}\left(\Omega\right)\subset L^{q}\left(\Omega\right)$.\\
$\hspace*{1em}$For the treatment of parabolic equations, we now give
the following
\end{singlespace}
\begin{defn}
\begin{singlespace}
\noindent A map $g\in L^{q}\left(\Omega\times(t_{0},t_{1})\right)$
belongs to the space $L^{q}\left(t_{0},t_{1};B_{q,\infty}^{\alpha}\left(\Omega\right)\right)$
if and only if
\[
\int_{t_{0}}^{t_{1}}\left(\sup_{h\in\mathbb{R}^{n}}\int_{\Omega}\frac{\left|g(x+h,t)-g(x,t)\right|^{q}}{\left|h\right|^{\alpha q}}\cdot\mathbf{{1}}_{\Omega}(x+h)\,dx\right)\,dt<\infty.
\]
\end{singlespace}
\end{defn}

\begin{singlespace}
\noindent $\hspace*{1em}$In this paper, we use the corresponding
local version of this space, which is denoted by the subscript ``\textit{loc}''.
Actually, in this framework we write $g\in L_{loc}^{q}\left(0,T;B_{q,\infty,loc}^{\alpha}\left(\Omega\right)\right)$
if and only if $g\in L^{q}\left(t_{0},t_{1};B_{q,\infty}^{\alpha}\left(\Omega'\right)\right)$
for all domains $\Omega'\times(t_{0},t_{1})\Subset\Omega_{T}$. Furthermore,
we shall also use the following notation, which is typical of Bochner
spaces: 
\[
\Vert g\Vert_{L^{q}\left(t_{0},\,t_{1}\,;\,B_{q,\infty}^{\alpha}\left(\Omega'\right)\right)}:=\left(\int_{t_{0}}^{t_{1}}\Vert g(\cdot,t)\Vert_{B_{q,\infty}^{\alpha}\left(\Omega'\right)}^{q}\,dt\right)^{\frac{1}{q}}.
\]

\end{singlespace}
\begin{singlespace}

\section{Proofs of the results\label{sec:proofs}}
\end{singlespace}

\begin{singlespace}
\noindent $\hspace*{1em}$We will now use the well-known difference
quotients method in the spatial directions, as well as the properties
of the function $H_{\lambda}$ and Proposition \ref{prop:Nikolskii}
to establish the\\

\noindent \begin{proof}[\bfseries{Proof of Theorem~\ref{thm:mainth}}]
By a slight abuse of notation, for $g\in L_{loc}^{1}\left(\Omega_{T},\mathbb{R}^{N}\right)$
and $i\in\left\{ 1,\ldots,n\right\} $, $h\neq0$, we set (when $x+he_{i}\in\Omega$)
\[
\tau_{h}g(x,t)\equiv\tau_{i,h}g(x,t):=g(x+he_{i},t)-g(x,t),
\]
\[
\Delta_{h}g(x,t)\equiv\Delta_{i,h}g(x,t):=\frac{g(x+he_{i},t)-g(x,t)}{h},
\]
where $e_{i}$ is the unit vector in the direction $x_{i}$. \\
$\hspace*{1em}$Since $u$ is a weak solution of equation (\ref{eq:1}),
we have
\[
\int_{\Omega_{T}}\left(u\cdot\partial_{t}\varphi-\langle H_{p-1}(Du),D\varphi\rangle\right)\,dz\,=\,-\int_{\Omega_{T}}f\varphi\,dz,
\]
for every test function $\varphi\in C_{0}^{\infty}(\Omega_{T})$.
Replacing $\varphi$ by $\tau_{-h}\varphi$, where $0<\left|h\right|<\mathrm{dist}(\mathrm{supp}\,\varphi,\partial\Omega_{T})$,
by virtue of the properties of the finite difference operator, we
get\smallskip{}
\[
\int_{\Omega_{T}}\left(\tau_{h}u\cdot\partial_{t}\varphi-\langle\tau_{h}H_{p-1}(Du),D\varphi\rangle\right)\,dz\,=\,-\int_{\Omega_{T}}\tau_{h}f\cdot\varphi\,dz.
\]
We now replace $\varphi$ by $\varphi_{\varepsilon}\equiv\phi_{\varepsilon}\ast\varphi$
in the previous equation, where $\left\{ \phi_{\varepsilon}\right\} $,
$\varepsilon>0$, denotes the family of standard, non-negative, radially
symmetric mollifiers in $\mathbb{R}^{n+1}$. This yields, for $0<\varepsilon\ll1$\smallskip{}
\[
\int_{\Omega_{T}}\left((\tau_{h}u)_{\varepsilon}\cdot\partial_{t}\varphi-\langle(\tau_{h}H_{p-1}(Du))_{\varepsilon},D\varphi\rangle\right)\,dz\,=\,-\int_{\Omega_{T}}(\tau_{h}f)_{\varepsilon}\cdot\varphi\,dz.
\]
Now, in the last equation we choose the test function $\varphi\equiv\Phi(\tau_{h}u)_{\varepsilon}$,
where $\Phi\in C_{0}^{\infty}(\Omega_{T})$ is a cut-off function
which will be specified later. After an integration by parts and then
letting $\varepsilon\searrow0$, we obtain\begin{equation}\label{eq:est0}
\begin{split}
&-\,\frac{1}{2}\int_{\Omega_{T}}\left|\tau_{h}u\right|^{2}\partial_{t}\Phi\,\,dz\,+\int_{\Omega_{T}}\Phi\langle\tau_{h}H_{p-1}(Du),D\tau_{h}u\rangle\,dz\\
&=\,-\int_{\Omega_{T}}\langle\tau_{h}H_{p-1}(Du),D\Phi\rangle\tau_{h}u\,\,dz\,+\int_{\Omega_{T}}\tau_{h}f\cdot\Phi\cdot\tau_{h}u\,\,dz.
\end{split}
\end{equation}

\noindent Note that an approximation argument yields the same identity
for any $\Phi\in W^{1,\infty}(\Omega_{T})$ with compact support in
$\Omega_{T}$ and any sufficiently small $h\in\mathbb{R}\setminus\left\{ 0\right\} $.
In what follows, we will denote by $c_{k}$ some positive constants
which do not depend on $h$.\\
$\hspace*{1em}$Now, let us consider a parabolic cylinder $Q_{\varrho}(z_{0})\subset Q_{2\varrho}(z_{0})\Subset\Omega_{T}$.
For a fixed time $t_{1}\in(t_{0}-\varrho^{2},t_{0})$ and $\delta\in(0,t_{0}-t_{1})$,
we choose $\Phi(x,t)=\widetilde{\chi}(t)\chi(t)\eta^{2}(x)$ with
$\chi\in W^{1,\infty}\left((0,T),\left[0,1\right]\right)$, $\chi\equiv0$
on $(0,t_{0}-\varrho^{2})$ and $\partial_{t}\chi\geq0$, $\eta\in C_{0}^{\infty}\left(B_{\varrho}(x_{0}),\left[0,1\right]\right)$,
and with the Lipschitz continuous function $\widetilde{\chi}:(0,T)\rightarrow\mathbb{R}$
defined by
\[
\widetilde{\chi}(t)=\begin{cases}
\begin{array}{c}
1\\
\mathrm{affine}\\
0
\end{array} & \begin{array}{c}
\mathrm{if}\,\,t\leq t_{1},\quad\quad\quad\;\,\,\\
\mathrm{if}\,\,t_{1}<t<t_{1}+\delta,\\
\mathrm{if}\,\,t\geq t_{1}+\delta.\quad\quad
\end{array}\end{cases}
\]

\noindent With such a choice of $\Phi$, equation \eqref{eq:est0}
turns into\begin{align*}
&-\,\frac{1}{2}\int_{\Omega_{T}}\left|\tau_{h}u\right|^{2}\eta^{2}(x)\chi(t)\partial_{t}\widetilde{\chi}(t)\,\,dz\,-\,\frac{1}{2}\int_{\Omega_{T}}\left|\tau_{h}u\right|^{2}\eta^{2}(x)\widetilde{\chi}(t)\partial_{t}\chi(t)\,\,dz\\
&+\int_{\Omega_{T}}\widetilde{\chi}(t)\chi(t)\eta^{2}(x)\langle\tau_{h}H_{p-1}(Du),D\tau_{h}u\rangle\,dz\\
&\,\,\,\,=-2\int_{\Omega_{T}}\widetilde{\chi}(t)\chi(t)\eta(x)\langle\tau_{h}H_{p-1}(Du),D\eta\rangle\tau_{h}u\,\,dz\,+\int_{\Omega_{T}}(\tau_{h}f)(\tau_{h}u)\widetilde{\chi}(t)\chi(t)\eta^{2}(x)\,dz.
\end{align*}Letting $\delta\rightarrow0$ in the previous equation, we get\\
\begin{equation}\label{eq:est1}
\begin{split}
\frac{1}{2}&\int_{B_{\varrho}(x_{0})}\chi(t_{1})\eta^{2}(x)\left|\tau_{h}u(x,t_{1})\right|^{2}dx\,+\int_{Q^{t_{1}}}\chi(t)\eta^{2}(x)\langle\tau_{h}H_{p-1}(Du),D\tau_{h}u\rangle\,dz\\
&=-2\int_{Q^{t_{1}}}\chi(t)\eta(x)\langle\tau_{h}H_{p-1}(Du),D\eta\rangle\tau_{h}u\,\,dz\,+\int_{Q^{t_{1}}}(\tau_{h}f)(\tau_{h}u)\chi(t)\eta^{2}(x)\,dz\\
&\,\,\,\,\,\,\,+\,\frac{1}{2}\int_{Q^{t_{1}}}(\partial_{t}\chi)\eta^{2}\left|\tau_{h}u\right|^{2}dz,
\end{split}
\end{equation}\\
for every $h\in\mathbb{R}\setminus\left\{ 0\right\} $ such that $\left|h\right|<\varrho/4$
and for almost every $t_{1}\in(t_{0}-\varrho^{2},t_{0})$, where we
have used the abbreviation $Q^{t_{1}}=B_{\varrho}(x_{0})\times(t_{0}-\varrho^{2},t_{1})$.\\
Now, by Lemma \ref{lem:Brasco} we have\\
\begin{equation}
\frac{4}{p^{2}}\int_{Q^{t_{1}}}\chi(t)\eta^{2}(x)\left|\tau_{h}H_{\frac{p}{2}}(Du)\right|^{2}dz\,\leq\,\int_{Q^{t_{1}}}\chi(t)\eta^{2}(x)\langle\tau_{h}H_{p-1}(Du),D\tau_{h}u\rangle\,dz,\label{eq:est2}
\end{equation}
and\begin{align*}
&\left|2\int_{Q^{t_{1}}}\chi(t)\eta(x)\langle\tau_{h}H_{p-1}(Du),D\eta\rangle\tau_{h}u\,\,dz\right|\\
&\leq2(p-1)\int_{Q^{t_{1}}}\chi(t)\eta(x)\left(\left|H_{\frac{p}{2}}(Du(x+he_{i},t))\right|^{\frac{p-2}{p}}+\left|H_{\frac{p}{2}}(Du)\right|^{\frac{p-2}{p}}\right)\left|\tau_{h}H_{\frac{p}{2}}(Du)\right|\left|D\eta\right|\left|\tau_{h}u\right|\,dz.
\end{align*}Using Young's inequality with exponents $\left(2,2\right)$ in the
right-hand side of the previous estimate, we obtain\begin{equation}\label{eq:est3}
\begin{split}
&\left|2\int_{Q^{t_{1}}}\chi(t)\eta(x)\langle\tau_{h}H_{p-1}(Du),D\eta\rangle\tau_{h}u\,\,dz\right|\\
&\leq\,\frac{2(p-1)^{2}}{\sigma}\int_{Q^{t_{1}}}\chi(t)\left|D\eta\right|^{2}\left(\left|H_{\frac{p}{2}}(Du(x+he_{i},t))\right|^{\frac{p-2}{p}}+\left|H_{\frac{p}{2}}(Du)\right|^{\frac{p-2}{p}}\right)^{2}\left|\tau_{h}u\right|^{2}dz\\
&\,\,\,\,\,\,\,+\,\frac{\sigma}{2}\int_{Q^{t_{1}}}\chi(t)\eta^{2}(x)\left|\tau_{h}H_{\frac{p}{2}}(Du)\right|^{2}dz,
\end{split}
\end{equation}where $\sigma>0$ will be chosen later. Moreover, by Hölder's inequality
with exponents $p$ and $p'$, as well as by the properties of the
difference quotients, of $\chi$, $\eta$, $u$ and $f$, we find
that\begin{equation}\label{eq:est4}
\begin{split}
&\left|\int_{Q^{t_{1}}}(\tau_{h}f)(\tau_{h}u)\chi(t)\eta^{2}(x)\,dz\right|\\
&\,\,\,\,\,\leq\,c_{1}(n)\,\left|h\right|\,\Vert Du\Vert_{L^{p}(Q_{2\varrho})}\left(\int_{Q^{t_{1}}}\left|\tau_{h}f\right|^{p'}dz\right)^{\frac{1}{p'}}=\,c_{1}\,\left|h\right|^{\alpha+1}\,\Vert Du\Vert_{L^{p}(Q_{2\varrho})}\left(\int_{Q^{t_{1}}}\frac{\left|\tau_{h}f\right|^{p'}}{\left|h\right|^{\alpha p'}}\,dz\right)^{\frac{1}{p'}}\\
&\,\,\,\,\,\leq c_{1}\,\left|h\right|^{\alpha+1}\,\Vert Du\Vert_{L^{p}(Q_{2\varrho})}\left[\int_{t_{0}\,-\,4\varrho^{2}}^{t_{0}}\left(\sup_{\left|y\right|<\frac{\varrho}{4}}\int_{B_{2\varrho}(x_{0})}\frac{\left|f(x+y,t)-f(x,t)\right|^{p'}}{\left|y\right|^{\alpha p'}}\cdot\mathbf{{1}}_{B_{2\varrho}(x_{0})}(x+y)\,dx\right)dt\right]^{\frac{1}{p'}}\\
&\,\,\,\,\,\leq\,c_{1}\,\left|h\right|^{\alpha+1}\,\Vert Du\Vert_{L^{p}(Q_{2\varrho})}\left(\int_{t_{0}\,-\,4\varrho^{2}}^{t_{0}}\Vert f(\cdot,t)\Vert_{B_{p',\infty}^{\alpha}(B_{2\varrho})}^{p'}\,dt\right)^{\frac{1}{p'}}\\
&\,\,\,\,\,=\,c_{1}\,\left|h\right|^{\alpha+1}\,\Vert Du\Vert_{L^{p}(Q_{2\varrho})}\,\Vert f\Vert_{L^{p'}\left(t_{0}\,-\,4\varrho^{2},\,t_{0}\,;\,B_{p',\infty}^{\alpha}\left(B_{2\varrho}(x_{0})\right)\right)},
\end{split}
\end{equation}for every $h\in\mathbb{R}\setminus\left\{ 0\right\} $ such that $\left|h\right|<\varrho/4$.
Joining estimates \eqref{eq:est1}, (\ref{eq:est2}), \eqref{eq:est3}
and \eqref{eq:est4}, and choosing $\sigma=4/p^{2}$, we arrive at\begin{align*}
&\int_{B_{\varrho}(x_{0})}\chi(t_{1})\eta^{2}(x)\left|\tau_{h}u(x,t_{1})\right|^{2}dx\,+\int_{Q^{t_{1}}}\chi(t)\eta^{2}(x)\left|\tau_{h}H_{\frac{p}{2}}(Du)\right|^{2}dz\\
&\leq c_{2}(p)\int_{Q^{t_{1}}}\left[\chi(t)\left|D\eta\right|^{2}\left(\left|H_{\frac{p}{2}}(Du(x+he_{i},t))\right|^{\frac{p-2}{p}}+\left|H_{\frac{p}{2}}(Du)\right|^{\frac{p-2}{p}}\right)^{2}+(\partial_{t}\chi)\eta^{2}\right]\left|\tau_{h}u\right|^{2}dz\\
&\,\,\,\,\,\,\,+\,c_{3}(n,p)\,\left|h\right|^{\alpha+1}\,\Vert Du\Vert_{L^{p}(Q_{2\varrho})}\,\Vert f\Vert_{L^{p'}\left(t_{0}\,-\,4\varrho^{2},\,t_{0}\,;\,B_{p',\infty}^{\alpha}\left(B_{2\varrho}(x_{0})\right)\right)},
\end{align*}which holds for almost every $t_{1}\in(t_{0}-\varrho^{2},t_{0})$.\\
We now choose a cut-off function $\eta\in C_{0}^{\infty}\left(B_{\varrho}(x_{0})\right)$
with $\eta\equiv1$ on $B_{\varrho/2}(x_{0})$ such that $0\leq\eta\leq1$
and $\left|D\eta\right|\leq C/\varrho$. For the cut-off function
in time, we choose the piecewise affine function $\chi:(0,T)\rightarrow\left[0,1\right]$
with
\[
\chi\equiv0\,\,\,\,\,\,\mathrm{on}\,\,\,(0,t_{0}-\varrho^{2}),\,\,\,\,\,\,\,\,\,\,\chi\equiv1\,\,\,\,\,\,\mathrm{on}\,\,\,(t_{0}-(\varrho/2)^{2},T)
\]
and
\[
\partial_{t}\chi\equiv\frac{4}{3\varrho^{2}}\,\,\,\,\,\,\mathrm{on}\,\,\,(t_{0}-\varrho^{2},t_{0}-(\varrho/2)^{2}).
\]
Dividing both sides of the previous estimate by $\left|h\right|$
and using the properties of $\chi$ and $\eta$, we obtain\begin{equation}\label{eq:est5}
\begin{split}
&\sup_{t_{0}-(\varrho/2)^{2}<t<t_{0}}\int_{B_{\varrho/2}(x_{0})}\left|\frac{\tau_{i,h}u(x,t)}{\left|h\right|^{1/2}}\right|^{2}dx\,\,+\,\int_{Q_{\varrho/2}(z_{0})}\left|\frac{\tau_{i,h}H_{p/2}(Du)}{\left|h\right|^{1/2}}\right|^{2}dz\\
&\,\,\,\,\,\,\,\leq \,c_{4}(p)\,\left|h\right|\,\varrho^{-2}\int_{Q_{\varrho}(z_{0})}\left[\left(\left|H_{\frac{p}{2}}(Du(x+he_{i},t))\right|^{\frac{p-2}{p}}+\left|H_{\frac{p}{2}}(Du)\right|^{\frac{p-2}{p}}\right)^{2}+1\right]\left|\Delta_{h}u\right|^{2}dz\\
&\,\,\,\,\,\,\,\,\,\,\,\,\,\,+\,c_{3}(n,p)\,\left|h\right|^{\alpha}\,\Vert Du\Vert_{L^{p}(Q_{2\varrho})}\,\Vert f\Vert_{L^{p'}\left(t_{0}\,-\,4\varrho^{2},\,t_{0}\,;\,B_{p',\infty}^{\alpha}\left(B_{2\varrho}(x_{0})\right)\right)}.
\end{split}
\end{equation}Now we set
\begin{equation}
I:=\int_{Q_{\varrho}(z_{0})}\left[\left(\left|H_{\frac{p}{2}}(Du(x+he_{i},t))\right|^{\frac{p-2}{p}}+\left|H_{\frac{p}{2}}(Du)\right|^{\frac{p-2}{p}}\right)^{2}+1\right]\left|\Delta_{h}u\right|^{2}dz\label{eq:intI}
\end{equation}
and we will assume that $p>2$. By Hölder's inequality with exponents
$\left(\frac{p}{2},\frac{p}{p-2}\right)$, as well as by the properties
of the difference quotients, we can control $I$ as follows\begin{equation}\label{eq:est6}
\begin{split}
I\,&\leq\, c_{5}(n)\left(\int_{Q_{2\varrho}}\left|Du\right|^{p}dz\right)^{\frac{2}{p}}\left(\int_{Q_{2\varrho}}\left[\left|H_{\frac{p}{2}}(Du)\right|^{\frac{2(p-2)}{p}}+1\right]^{\frac{p}{p-2}}dz\right)^{\frac{p-2}{p}}\\
&\leq\, c_{6}(n,p)\left(\int_{Q_{2\varrho}}\left|Du\right|^{p}dz\right)^{\frac{2}{p}}\left(\int_{Q_{2\varrho}}\left[\left|Du\right|^{p}+1\right]\,dz\right)^{\frac{p-2}{p}}\\
&\leq\,c_{6}(n,p)\,\int_{Q_{2\varrho}}\left(\left|Du\right|^{p}+1\right)\,dz,
\end{split}
\end{equation}provided that $\left|h\right|<\varrho/4$. Combining estimates \eqref{eq:est5}
and \eqref{eq:est6}, we then have \begin{equation}\label{eq:est7}
\begin{split}
&\sup_{t_{0}-(\varrho/2)^{2}<t<t_{0}}\int_{B_{\varrho/2}(x_{0})}\left|\frac{\tau_{i,h}u(x,t)}{\left|h\right|^{(\alpha+1)/2}}\right|^{2}dx\,\,+\,\int_{Q_{\varrho/2}(z_{0})}\left|\frac{\tau_{i,h}H_{p/2}(Du)}{\left|h\right|^{(\alpha+1)/2}}\right|^{2}dz\\
&\,\,\,\,\,\,\,\leq\,c\,\left|h\right|^{1-\alpha}\,\varrho^{-2}\int_{Q_{2\varrho}}\left(\left|Du\right|^{p}+1\right)\,dz\,+\,c\,\,\Vert Du\Vert_{L^{p}(Q_{2\varrho})}\,\Vert f\Vert_{L^{p'}\left(t_{0}\,-\,4\varrho^{2},\,t_{0}\,;\,B_{p',\infty}^{\alpha}\left(B_{2\varrho}(x_{0})\right)\right)},
\end{split}
\end{equation}with $c\equiv c(n,p)>0$. Since the previous estimate holds for every
$h\in\mathbb{R}\setminus\left\{ 0\right\} $ such that $\left|h\right|<\varrho/4$,
from \eqref{eq:est7} we obtain\begin{align*}
\int_{Q_{\varrho/2}(z_{0})}\left|\frac{\tau_{i,h}H_{p/2}(Du)}{\left|h\right|^{(\alpha+1)/2}}\right|^{2}dz\,&\leq\,4^{\alpha-1}\,c\,\varrho^{-1-\alpha}\int_{Q_{2\varrho}}\left(\left|Du\right|^{p}+1\right)\,dz\\
&\,\,\,\,\,\,\,+\,c\,\,\Vert Du\Vert_{L^{p}(Q_{2\varrho})}\,\Vert f\Vert_{L^{p'}\left(t_{0}\,-\,4\varrho^{2},\,t_{0}\,;\,B_{p',\infty}^{\alpha}\left(B_{2\varrho}(x_{0})\right)\right)},
\end{align*}from which we can deduce 
\[
\int_{Q_{\varrho/2}(z_{0})}\left|\tau_{i,h}H_{\frac{p}{2}}(Du)\right|^{2}dz\,\leq\,M\left|h\right|^{\alpha+1},
\]
for some finite positive constant $M$ depending on $p$, $n$, $\varrho$,
$\alpha$, $\Vert f\Vert_{L^{p'}\left(t_{0}\,-\,4\varrho^{2},\,t_{0}\,;\,B_{p',\infty}^{\alpha}(B_{2\varrho})\right)}$
and $\Vert Du\Vert_{L^{p}(\Omega_{T})}$, but not on $h$. Note that
the above estimate holds for every $i\in\left\{ 1,\ldots,n\right\} $
and every sufficiently small $h\in\mathbb{R}\setminus\left\{ 0\right\} $.
Therefore, using Proposition \ref{prop:Nikolskii} with the choices
$G=H_{p/2}(Du)$, $q=2$ and $\theta=\frac{\alpha+1}{2}$, as well
as a standard covering argument, we infer that 
\[
H_{\frac{p}{2}}(Du)\,\in\,L_{loc}^{2}\left(0,T;W_{loc}^{s,2}\left(\Omega,\mathbb{R}^{n}\right)\right)\,\,\,\,\,\,\,\,\,\mathrm{for\,\,all}\,\,s\in\left(0,\frac{\alpha+1}{2}\right).
\]
$\hspace*{1em}$Finally, when $p=2$, arguing in a similar fashion
we reach the same conclusions. This completes the proof.\end{proof}

\noindent $\hspace*{1em}$We are now in a position to give the\\

\noindent \begin{proof}[\bfseries{Proof of Corollary~\ref{cor:cor1}}]
By virtue of Theorems \ref{thm:mainth}, \ref{thm:Lipschitz} and
\ref{thm:fracemb}, for any fixed $s\in\left(0,\frac{\alpha+1}{2}\right)$
we have
\[
H_{\frac{p}{2}}(Du)\,\in\,L^{2}\left(t_{0},t_{1};L^{\gamma}\left(\Omega',\mathbb{R}^{n}\right)\right)\,\,\,\,\,\,\,\,\,\mathrm{for\,\,every}\,\,\gamma\in\left[1,\frac{2n}{n-2s}\right],
\]
for any open Lipschitz set $\Omega'\Subset\Omega$ and any $(t_{0},t_{1})\subset(0,T)$.
Now, let us observe that the function $g_{n}:\left(0,\frac{\alpha+1}{2}\right)\rightarrow\mathbb{R}$
defined by
\[
g_{n}(s)=\frac{2n}{n-2s}
\]
is continuous and strictly increasing, and $g_{n}(s)\nearrow\frac{2n}{n\,-\,\alpha\,-\,1}$
as $s\nearrow\frac{\alpha+1}{2}$. Hence, letting $s$ tend to $\frac{\alpha+1}{2}$
from below, we obtain that
\[
H_{\frac{p}{2}}(Du)\,\in\,L_{loc}^{2}\left(0,T;L_{loc}^{r}\left(\Omega,\mathbb{R}^{n}\right)\right)\,\,\,\,\,\,\,\,\,\mathrm{for\,\,all}\,\,r\in\left[1,\frac{2n}{n-\alpha-1}\right),
\]
since every open ball $B_{\varrho}\Subset\Omega$ is a set of class
$C^{0,1}$ with bounded boundary. In particular, for all Lipschitz
open subsets $\Omega'\Subset\Omega$ and all $(t_{0},t_{1})\subset(0,T)$
we have
\begin{equation}
\int_{t_{0}}^{t_{1}}\left(\int_{\Omega'}(\vert Du(x,t)\vert-1)_{+}^{rp/2}\,dx\right)^{\frac{2}{r}}dt\,=\,\int_{t_{0}}^{t_{1}}\Vert H_{\frac{p}{2}}(Du(\cdot,t))\Vert_{L^{r}(\Omega')}^{2}\,dt\,<+\infty\label{eq:highsumm}
\end{equation}
for any $r\in\left[1,\frac{2n}{n\,-\,\alpha\,-\,1}\right)$. Now notice
that 
\[
1\leq\,r\,<\frac{2n}{n-\alpha-1}\,\,\,\,\,\,\Longrightarrow\,\,\,\,\,\,1\leq\,\frac{p}{2}\,\leq\,\frac{rp}{2}\,<\,\frac{np}{n-\alpha-1}\,,
\]
since $p\geq2$. Therefore, from (\ref{eq:highsumm}) it follows that
\[
\int_{t_{0}}^{t_{1}}\Vert(\vert Du(\cdot,t)\vert-1)_{+}\Vert_{L^{q}(\Omega')}^{p}\,dt\,=\,\int_{t_{0}}^{t_{1}}\left(\int_{\Omega'}(\vert Du(x,t)\vert-1)_{+}^{q}\,dx\right)^{\frac{p}{q}}dt\,<+\infty
\]
for all $q\in\left[1,\frac{np}{n\,-\,\alpha\,-\,1}\right)$, for any
open Lipschitz set $\Omega'\Subset\Omega$ and any $(t_{0},t_{1})\subset(0,T)$.
This is sufficient to ensure that
\[
Du\,\in\,L_{loc}^{p}\left(0,T;L_{loc}^{q}\left(\Omega,\mathbb{R}^{n}\right)\right)\,\,\,\,\,\,\,\,\,\mathrm{for\,\,all}\,\,q\in\left[1,\frac{np}{n-\alpha-1}\right).
\]
\end{proof}

\noindent \smallskip{}

\noindent $\hspace*{1em}$\textbf{Acknowledgements. }The author is
a member of the \textit{Gruppo Nazionale per l'Analisi Matematica,
la Probabilità e le loro Applicazioni }(GNAMPA) of the \textit{Istituto
Nazionale di Alta Matematica} (INdAM).
\end{singlespace}

\begin{singlespace}

\lyxaddress{\noindent \textbf{$\quad$}\\
$\hspace*{1em}$\textbf{Pasquale Ambrosio}\\
Dipartimento di Matematica e Applicazioni ``R. Caccioppoli''\\
Università degli Studi di Napoli ``Federico II''\\
Via Cintia, 80126 Napoli, Italy.\\
\textit{E-mail address}: pasquale.ambrosio2@unina.it}
\end{singlespace}

\end{document}